%% file: hkmp10.tex
\documentclass[11pt,leqno]{article}
\usepackage{amsmath}
\usepackage{epsfig}
\usepackage{epic,eepic}
\def\disp{\displaystyle}

\def\Limsup{\mathop{{\rm Lim}\,{\rm sup}}}

\def\tto{\;{\lower 1pt \hbox{$\rightarrow$}}\kern -10pt
\hbox{\raise 2pt \hbox{$\rightarrow$}}\;}
\def\Hat{\widehat}

\def\Tilde{\widetilde}

\def\ra{\rangle}
\def\la{\langle}
\def\ve{\varepsilon}

\def\h{\hfill\Box}
\def\R{I\!\!R}
\def\N{I\!\!N}
\def\ox{\bar{x}}
\def\oy{\bar{y}}
\def\oz{\bar{z}}

\def\gph{\mbox{\rm gph}\,}
\def\epi{\mbox{\rm epi}\,}

\def\ker{\mbox{\rm ker}\,}

\def\h{\hfill\triangle}
\def\dn{\downarrow}
\def\O{\Omega}

\def\emp{\emptyset}
\def\st{\stackrel}

\def\lm{\lambda}
\def\gg{\gamma}

\def\dd{\delta}
\def\al{\alpha}

\def\N{I\!\!N}

\def\tilde{\widetilde}
\newcounter{lk}

\begin{document}
\begin{center}
\vspace*{0.3in} {\bf GENERALIZED NEWTON'S METHOD\\ BASED ON
GRAPHICAL DERIVATIVES}\\[3ex]
T. HOHEISEL\footnote{Institute of Applied Mathematics and
Statistics, University of W\"urzburg, Am Hubland, 97074 W\"urzburg, Germany
(hoheisel@mathematik.uni-wuerzburg.de,
kanzow@mathematik.uni-wuerzburg.de).}, C. KANZOW\footnotemark[1],
B. S. MORDUKHOVICH\footnote{Department of Mathematics, Wayne State
University, Detroit, MI 48202, USA (boris@math.wayne.edu,
pmhung@wayne.edu). The research of these authors was partially
supported by the US National Science Foundation under grants
DMS-0603846 and DMS-1007132 and by the Australian Research Council
under grant DP-12092508.} and H. PHAN\footnotemark[2]
\end{center}

\small{\bf Abstract.} This paper concerns developing a numerical
method of the Newton type to solve systems of nonlinear equations
described by nonsmooth continuous functions. We propose and
justify a new generalized Newton algorithm based on graphical
derivatives, which have never been used to derive a Newton-type
method for solving nonsmooth equations. Based on advanced
techniques of variational analysis and generalized
differentiation, we establish the well-posedness of the algorithm,
its local superlinear convergence, and its global convergence of
the Kantorovich type. Our convergence results hold with no
semismoothness assumption, which is illustrated by examples. The
algorithm and main results obtained in the paper are compared with
well-recognized semismooth and $B$-differentiable versions of
Newton's method for nonsmooth Lipschitzian equations.
\vspace*{0.05in}

{\bf Key words.} nonsmooth equations, optimization and variational
analysis, Newton's method, graphical derivatives and
coderivatives, local and global convergence\vspace*{0.05in}

{\bf AMS subject classification.} 49J53, 65K15, 90C30

\newtheorem{Theorem}{Theorem}[section]
\newtheorem{Proposition}[Theorem]{Proposition}
\newtheorem{Remark}[Theorem]{Remark}
\newtheorem{Lemma}[Theorem]{Lemma}
\newtheorem{Corollary}[Theorem]{Corollary}
\newtheorem{Definition}[Theorem]{Definition}
\newtheorem{Example}[Theorem]{Example}
\newtheorem{Algorithm}[Theorem]{Algorithm}
\renewcommand{\theequation}{\thesection.\arabic{equation}}
\normalsize

\section{Introduction}
\setcounter{equation}{0}

Newton's method is one of the most powerful and useful methods in
optimization and in the related area of solving systems of
nonlinear equations
\begin{equation}\label{eq:ne}
H(x)=0
\end{equation}
defined by continuous vector-valued mappings $H\colon\R^n\to\R^n$.
In the classical setting when $H$ is a continuously differentiable
(smooth, $C^1$) mapping, Newton's method builds the following
iteration procedure
\begin{equation}\label{eq:iteration}
x^{k+1}:=x^k+d^k\;\mbox{ for all }\;k=0,1,2,\ldots,
\end{equation}
where $x^0\in\R^n$ is a given starting point, and where
$d^k\in\R^n$ is a solution to the linear system of equations
(often called ``Newton equation")
\begin{equation}\label{eq:linearized}
H'(x^k)d=-H(x^k).
\end{equation}
A detailed analysis and numerous applications of the classical
Newton's method \eqref{eq:iteration}, \eqref{eq:linearized} and
its modifications can be found, e.g., in the books \cite{DeS
83,Kel 03,OrR 70} and the references therein.

However, in the vast majority of applications---including those to
optimization, variational inequalities, complementarity and
equilibrium problems, etc.---the underlying mapping $H$ in
\eqref{eq:ne} is nonsmooth. Indeed, the aforementioned
optimization-related models and their extensions can be written
via Robinson's formalism of ``generalized equations," which in
turn can be reduced to standard equations of the form above
(using, e.g., the projection operator) while with {\em
intrinsically nonsmooth} mappings $H$; see \cite{FaP 03,LPR
96,rob91,PaQ 93} for more details, discussions, and references.

Robinson originally proposed (see \cite{Rob 87} and also \cite{Rob 94}
based on his earlier preprint) a point-based approximation approach
to solve nonsmooth equations \eqref{eq:ne}, which then was developed
by his student Josephy \cite{Jos 79} to extend Newton's method for
solving variational inequalities and complementarity problems. Other
approaches replace the classical derivative $H'(x_k)$ in the Newton
equation \eqref{eq:linearized} by some generalized derivatives. In
particular, the $B$-differentiable Newton method developed by Pang
\cite{Pan 90,Pan 91} uses the iteration scheme \eqref{eq:iteration}
with $d^k$ being a solution to the subproblem
\begin{equation}\label{b1}
H'(x^k; d)=-H(x^k),
\end{equation}
where $H'(x^k;d)$ denotes the classical directional derivative of
$H$ at $x^k$ in the direction $d$. Besides the existence of the
classical directional derivative in \eqref{b1}, a number of strong
assumptions are imposed in \cite{Pan 90, Pan 91} to establish
appropriate convergence results; see Section~5 below for more
discussions and comparisons.

In another approach developed by Kummer \cite{k88} and Qi and Sun
\cite{QiS 93}, the direction $d^k$ in \eqref{eq:iteration} is taken
as a solution to the linear system of equations
\begin{equation}\label{c1}
A_k d=-H(x^k),
\end{equation}
where $A_k$ is an element of Clarke's generalized Jacobian
$\partial_C H(x_k)$ of a Lipschitz continuous mapping $H$. In
\cite{Qi 93}, Qi suggested to replace $A_k\in\partial_C H(x^k)$ in
\eqref{c1} by the choice of $A_k$ from the so-called
$B$-subdifferential $\partial_B H(x^k)$ of $H$ at $x^k$, which is
a proper subset of $\partial_C H(x^k)$; see Section~4 for more
details. We also refer the reader to \cite{FaP 03,kk, Rob 94} and
bibliographies therein for wide overviews, historical remarks, and
other developments on Newton's method for nonsmooth Lipschitz
equations as in \eqref{eq:ne} and to \cite{kff} for some recent
applications.

It is proved in \cite{QiS 93} and \cite{Qi 93} that the Newton
type method based on implementing the generalized Jacobian and
$B$-subdifferential in \eqref{c1}, respectively, superlinearly
converges to a solution of \eqref{eq:ne} for a class of {\em
semismooth} mappings $H$; see Section~4 for the definition and
discussions. This subclass of Lipschitz continuous and
directionally differentiable mappings is rather broad and useful
in applications to optimization-related problems. However, not
every mapping arising in applications (from both theoretical and
practical viewpoints) is either directionally differentiable or
Lipschitz continuous. The reader can find valuable classes of
functions and mappings of this type in \cite{mor06,rw} and
overwhelmingly in spectral function analysis, eigenvalue
optimization, studying of roots of polynomials, stability of
control systems, etc.; see, e.g., \cite{blo} and the references
therein.\vspace*{0.05in}

The main goal and achievements of this paper are as follows. We
propose a new Newton-type algorithm to solve nonsmooth equations
\eqref{eq:ne} described by general continuous mappings $H$ that is
based on {\em graphical derivatives}. It reduces to the classical
Newton method \eqref{eq:linearized} when $H$ is smooth, being
different from previously known versions of Newton's method in the
case of Lipschitz continuous mappings $H$. Based on advanced tools
of variational analysis involving {\em metric regularity} and {\em
coderivatives}, we justify well-posedness of the new algorithm and
its {\em superlinear local} and {\em global} (of the Kantorovich
type) convergence under verifiable assumptions that hold for
semismooth mappings but are {\em not} restricted to them. Detailed
comparisons of our algorithm and results with the semismooth and
$B$-differentiable Newton methods are given and certain
improvements of these methods are justified.

Note metric regularity and related concepts of variational
analysis has been employed in the analysis and justification of
numerical algorithms starting with Robinson's seminal
contribution; see, e.g., \cite{ag,llm,mpr} and their references
for the recent account. However, we are not familiar with any
usage of graphical derivatives and coderivatives for these
purposes.\vspace*{0.05in}

The rest of the paper is organized as follows. In Section~2 we
present basic definitions and preliminaries from variational
analysis and generalized differentiation widely used for
formulations and proofs of the main results.

Section~3 is devoted to the description of the new generalized
Newton algorithm with justifying its well-posedness/solvability
and establishing its superlinear local and global convergence
under appropriate assumptions on the underlying mapping $H$.

In Section~4 we compare our algorithm with the scheme of
\eqref{c1}. We also discuss in detail the major assumptions made
in Section~3 deriving sufficient conditions for their fulfillment
and comparing them with those in the semismooth Newton methods.

Section~5 contains applications of our algorithm to the
$B$-differentiable Newton method \eqref{b1} with largely relaxed
assumptions in comparison with known ones. In Section~6 we give
some concluding remarks and discussions on further research.
\vspace*{0.05in}

Our notation is basically standard in variational analysis and
numerical optimization; cf. \cite{FaP 03,mor06,rw}. Recall that,
given a set-valued mapping $F\colon\R^n\tto\R^m$, the expression
\begin{eqnarray}\label{pk}
\begin{array}{ll}
\disp\Limsup_{x\to\ox}F(x):=\Big\{y\in\R^m\Big|&\exists\,x_k\to\ox\;\mbox{
and }\;y_k\to y\;\mbox{ as }\;k\to\infty\;\mbox{ with}\\
&y_k\in F(x_k)\;\mbox{ for all }\;k\in\N:=\{1,2,\ldots\}\Big\}
\end{array}
\end{eqnarray}
defines the {\em Painlev\'e-Kuratowski upper/outer limit} of $F$
as $x\to\ox$. Let us also mention that the symbols ${\rm
cone}\,\O$ and ${\rm co}\,\O$ stand, respectively, for the conic
hull and convex hull of the set in question, that ${\rm
dist}(x;\O)$ denotes the Euclidean distance between a point
$x\in\R^n$ and a set $\O$, and that the notation $A^T$ signifies
the matrix transposition. As usual, $B_\ve(\ox)$ stands for the
closed ball centered at $\ox$ with radius $\ve>0$.

\section{Tools of Variational Analysis}
\setcounter{equation}{0}

In this section we briefly review some constructions and results
from variational analysis and generalized differentiation widely
used in what follows. The reader may consult the texts
\cite{AuF 90,mor06,rw,s} for more details and additional material.

Given a nonempty set $\O\subset\R^n$ and a point $\ox\in\O$, the
(Bouligand-Severi) {\em tangent/contingent cone} to $\O$ at $\ox$ is
defined by
\begin{equation}\label{tan}
   T(\ox;\O):=\disp\Limsup_{t\dn 0}\frac{\O-\ox}{t}
\end{equation}
via the outer limit \eqref{pk}. This cone is often nonconvex while its
polar/dual cone
\begin{equation}\label{fn}
   \Hat N(\ox;\O):=\big\{p\in\R^n\big|\;\la p,u\ra\le 0\;
   \mbox{ for all }\;u\in T(\ox;\O)\big\}
\end{equation}
is always convex and can be intrinsically described by
\begin{equation*}
   \Hat N(\ox;\O)=\Big\{p\in\R^n\Big|\;\limsup_{x\st{\O}{\to}\ox}
   \frac{\la p,x-\ox \ra}{\|x-\ox\|} \le 0\Big\},\quad\ox\in\O,
\end{equation*}
where the symbol $x\st{\O}{\to}\ox$ signifies that $x\to\ox$ with
$x\in\O$. The construction \eqref{fn} is known as the {\em
prenormal cone} or the {\em Fr\'echet/regular normal cone} to $\O$
at $\ox\in\O$. For convenience we put $\Hat N(\ox;\O)=\emp$ if
$\ox\notin\O$. Observe that the prenormal cone \eqref{fn} may not
have natural properties of generalized normals in the case of
nonconvex sets $\O$; e.g., it often happens that $\Hat
N(\ox;\O)=\{0\}$ when $\ox$ is a boundary point of $\O$ and the
cone \eqref{fn} does not possesses required calculus rules. The
situation is dramatically improved when we consider a robust
regularization of \eqref{fn} via the outer limit \eqref{pk} and
arrive at the construction
\begin{equation}\label{nc}
   N(\ox;\O):=\Limsup_{x\to\ox}\Hat N(x;\O)
\end{equation}
known as the (limiting, basic, Mordukhovich) {\em normal cone} to
$\O$ at $\ox\in\O$. If $\O$ is locally closed around $\ox$, the
basic normal cone \eqref{nc} can be equivalently described as
\begin{equation*}
   N(\ox;\O)=\Limsup_{x\to\ox}\big[{\rm cone}\big(x-\Pi(x;\O)\big)\big],
   \quad\ox\in\O,
\end{equation*}
via the Euclidean projector $\Pi(\cdot;\O)$ on $\O$; this was in
fact the original definition of the normal cone in \cite{m76}.
Despite its nonconvexity, the normal cone \eqref{nc} and the
corresponding subdifferential and coderivative constructions for
extended-real-valued functions and set-valued mappings enjoy
comprehensive calculus rules, which are particularly based on
variational/extremal principles of variational analysis.
\vspace*{0.05in}

Consider next a set-valued mapping $F\colon\R^n\tto\R^m$ with the graph
\begin{equation*}
   \gph F:=\big\{(x,y)\in\R^n\times\R^m\big|\;y\in F(x)\big\}
\end{equation*}
and define the graphical derivative and coderivative constructions
generated by the tangent and normal cones, respectively. Given
$(\ox,\oy)\in\gph F$, the {\em graphical/contingent derivative} of $F$ at
$(\ox,\oy)$ is introduced in \cite{Aub 81} as a mapping
$DF(\ox,\oy)\colon\R^n\tto\R^m$ with the values
\begin{equation}\label{d}
   DF(\ox,\oy)(z):=\big\{w\in\R^m\big|\;(z,w)\in T\big((\ox,\oy);\gph
   F\big)\big\},\quad z\in\R^n,
\end{equation}
defined via the contingent cone \eqref{tan} to the graph of $F$ at
the point $(\ox,\oy)$; see \cite{AuF 90,rw} for various
properties, equivalent representation, and applications. The {\em
coderivative} of $F$ at $(\ox,\oy)\in\gph F$ is introduced in
\cite{mor80} as a mapping $D^*F(\ox,\oy)\colon\R^m\tto\R^n$ with
the values
\begin{equation}\label{cod}
   D^*F(\ox,\oy)(v):=\big\{u\in\R^n\big|\;(u,-v)\in
   N\big((\ox,\oy);\gph F\big)\big\},\quad v\in\R^m,
\end{equation}
defined via the normal cone \eqref{nc} to the graph of $F$ at
$(\ox,\oy)$; see \cite{mor06,rw} for extended calculus and a
variety of applications.  We drop $\oy$ in the graphical derivative and
coderivative notation when the mapping in question is
single-valued at $\ox$. Note that the graphical derivative and
coderivative constructions in \eqref{d} and \eqref{cod} are {\em
not dual} to each other, since the basic normal cone \eqref{nc} is
{\em nonconvex} and hence cannot be tangentially generated.

In this paper we employ, together with \eqref{d} and \eqref{cod},
the following modified derivative construction for mappings, which
seems to be new in generality although constructions of this
(radial, Dini-like) type have been widely used for
extended-real-valued functions.

\begin{Definition}\label{def:restrder} {\bf (restrictive graphical
derivative of mappings).}
Let $F\colon\R^n\tto\R^m$, and let $(\ox,\oy)\in\gph F$. Then a
set-valued mapping $\Tilde DF(\ox,\oy)\colon\R^n\tto\R^m$ given by
\begin{equation}\label{rd}
   \Tilde DF(\ox,\oy)(z):=\Limsup_{t\dn 0}\frac{F(\ox+tz)-\oy}{t},
   \quad z\in\R^n,
\end{equation}
is called the {\sc restrictive graphical derivative} of $F$ at $(\ox,\oy)$.
\end{Definition}

The next proposition collects some properties of the graphical
derivative \eqref{d} and its restrictive counterpart \eqref{rd}
needed in what follows.

\begin{Proposition}\label{prop:derivprop} {\bf (properties of graphical
derivatives).}
Let $F\colon\R^n\tto\R^m$, and let $(\ox,\oy)\in\gph F$. Then the
following assertions hold:

{\bf (a)} We have $\Tilde DF(\ox,\oy)(z)\subset DF(\ox,\oy)(z)$ for all
$z\in\R^n$.

{\bf (b)} There are inverse derivative relationships
\begin{equation*}
   DF(\ox,\oy)^{-1}=DF^{-1}(\oy,\ox)\;\mbox{ and }\;
   \Tilde DF(\ox,\oy)^{-1}=\Tilde DF^{-1}(\oy,\ox).
\end{equation*}

{\bf (c)} If $F$ is single-valued and locally Lipschitzian around $\ox$, then
\begin{equation*}
   \Tilde DF(\ox)(z)=DF(\ox)(z)\;\mbox{ for all }\;z\in\R^n.
\end{equation*}

{\bf (d)} If $F$ is single-valued and directionally differentiable at $\ox$,
then
\begin{equation*}
   \Tilde DF(\ox)(z)=\big\{F'(\ox;z)\big\}\;\mbox{ for all }\;z\in\R^n.
\end{equation*}

{\bf (e)} If $F$ is single-valued and G\^ateaux differentiable at $\ox$
with the G\^ateaux derivative $F'_G(\ox)$, then we have
\begin{equation*}
   \Tilde DF(\ox)(z)=\big\{F'_G(\ox)z\big\}\;\mbox{ for all }\;z\in\R^n.
\end{equation*}

{\bf (f)} If $F$ is single-valued and $($Fr\'echet$)$ differentiable at
$\ox$ with the derivative $F'(\ox)$, then
\begin{equation*}
   DF(\ox)(z)=\big\{F'(\ox)z\big\}\;\mbox{ for all }\;z\in\R^n.
\end{equation*}
\end{Proposition}

\noindent
{\bf Proof.}
It is shown in \cite[8(14)]{rw} that the graphical derivative \eqref{d}
admits the representation
\begin{equation}\label{gd}
   DF(\ox,\oy)(z) = \Limsup_{t\dn 0,\,h\to z}\frac{F(\ox+th)-\oy}{t},
   \quad z\in\R^n.
\end{equation}
The inclusion in (a) is an immediate consequence of
Definition~\ref{def:restrder} and representation \eqref{gd}.

The first equality in (b), observed from the very beginning \cite{Aub
81}, easily follows from definition \eqref{d}. We can similarly
check the second one in (b).

To justify the equality in (c), it remains to verify by (a) the
opposite inclusion `$\supset$' when $F$ is single-valued and
locally Lipschitzian around $\ox$. In this case fix $z\in\R^n$,
pick any $w\in DF(\ox)(z)$, and find by representation \eqref{gd}
sequences $h_k\to z$ and $t_k\dn 0$ such that
\begin{equation*}
   \frac{F(\ox+t_k h_k)-F(\ox)}{t_k}\to w\;\mbox{ as }\;k\to\infty.
\end{equation*}
The local Lipschitz continuity of $F$ around $\ox$ with constant $L\ge 0$
implies that
\begin{eqnarray*}
   \Big\|\frac{F(\ox+t_kh_k)-F(\ox)}{t_k}-\frac{F(\ox+t_k z)-F(\ox)}{t_k}\Big\|
   &=&\Big\|\frac{F(\ox+t_kh_k)-F(\ox+t_kz)}{t_k}\Big\|\\
   &\le& L\|h_k-z\Big\|
\end{eqnarray*}
for all $k\in\N$ sufficiently large, and hence we have the convergence
\begin{equation*}
   \frac{F(\ox+t_k z)-F(\ox)}{t_k}\to w\;\mbox{ as }\;k\to\infty.
\end{equation*}
Thus $w\in\Tilde DF(\ox)(z)$, which justifies (c). Assertions (d)
and (e) follow directly from the definitions. Finally, assertion
(f) is implied by (e) in the local Lipschitzian case (c) while it
can be easily derived from the (Fr\'echet) differentiability of
$F$ at $\ox$ with no Lipschitz assumption; see, e.g.,
\cite[Exercise~9.25(b)]{rw}. $\h$\vspace*{0.05in}

Proposition~\ref{prop:derivprop} reveals important differences
between the graphical derivative \eqref{d} and the coderivative
\eqref{cod}. Indeed, assertions (c) and (d) of this proposition
show that the graphical derivative of locally Lipschitzian and
{\em directionally differentiable} mappings $F\colon\R^n\to\R^m$
is always {\em single-valued}.  At the same time, the coderivative
single-valuedness for locally Lipschitzian mappings is equivalent
to the {\em strict/strong} Fr\'echet {\em differentiability} of
$F$ at the point in question; see \cite[Theorem~3.66]{mor06}. It
follows from the well-known formula
\begin{equation}\label{C2}
   {\rm co}D^*F(\ox)(z)=\big\{A^Tz\big|\;A\in\partial_C F(\ox)\big\}
\end{equation}
that the latter strict differentiability condition characterizes also the
single-valuedness of the generalized Jacobian of $F$ at $\ox$.

In fact, in the case of $F=(f_1,\ldots,f_m)\colon\R^n\to\R^m$
being locally Lipschitzian around $\ox$ the coderivative
\eqref{cod} admits the subdifferential description
\begin{equation}\label{Eq:scalarization}
   D^*F(\ox)(z)=\partial\Big(\sum_{i=1}^m z_i f_i\Big)(\ox)\;\mbox{
   for any }\;z=(z_1,\ldots,z_m)\in\R^m,
\end{equation}
where the (basic, limiting, Mordukhovich) {\em subdifferential}
$\partial f(\ox)$ of a general scalar function $f$ at $\ox$ is
defined geometrically by
\begin{equation}\label{sub}
   \partial f(\ox):=\big\{p\in\R^n\big|\;(p,-1)\in N\big((\ox,f(\ox));
   \epi f\big)\big\}
\end{equation}
via the normal cone \eqref{nc} to the epigraph $\epi
f:=\{(x,\mu)\in\R^{n+1}|\;\mu\ge f(x)\}$ and admits analytical
descriptions in terms of the outer limit \eqref{pk} of the
Fr\'echet/regular and proximal subdifferentials at points nearby;
see \cite{mor06,rw} with the references therein. Note also that
the basic subdifferential \eqref{sub} of a continuous function $f$
can be also described via the coderivative of $f$ by $\partial
f(\ox)=D^*f(\ox)(1)$; see
\cite[Theorem~1.80]{mor06}.\vspace*{0.05in}

Finally in this section, we recall the notion of metric regularity
and its coderivative characterization that play a significant role
in the paper. A mapping $F\colon\R^n\tto\R^m$ is {\em metrically
regular} around $(\ox,\oy)\in\gph F$ if there are neighborhoods
$U$ of $\ox$ and $V$ of $\oy$ as well as a number $\mu>0$ such
that
\begin{equation}\label{mr}
   {\rm dist}\big(x;F^{-1}(y)\big)\le\mu\,{\rm dist}\big(y;F(x)\big)\;
   \mbox{ for all }\;x\in U\;\mbox{ and }\;y\in V.
\end{equation}
Observe that it is sufficient to require the fulfillment of
\eqref{mr} just for those $y\in V$ satisfying the estimate
dist$(y;F(x))\le\gg$ for some $\gg>0$; see
\cite[Proposition~1.48]{mor06}.

We will see below that metric regularity is crucial for justifying
the well-posedness of our generalized Newton algorithm and
establishing its local and global convergence. It is also worth
mentioning that, in the opposite direction, a Newton-type method
(known as the Lyusternik-Graves iterative process) leads to
verifiable conditions for metric regularity of smooth mappings;
see, e.g., the proof of \cite[Theorem~1.57]{mor06} and the
commentaries therein. The latter procedure is replaced by
variational/extremal principles of variational analysis in the
case of nonsmooth and set-valued mappings under consideration;
cf.\ \cite{i,mor06,rw}.

In this paper we broadly use the following {\em coderivative
characterization} of the metric regularity property for an
arbitrary set-valued mapping $F$ with closed graph, known also as
the Mordukhovich criterion (see \cite[Theorem~3.6]{mor93},
\cite[Theorem~9.45]{rw}, and the references therein): $F$ is
metrically regular around $(\ox,\oy)$ {\em if and only if} the
inclusion
\begin{equation}\label{eq:coderivativecriterion}
0\in D^*F(\ox,\oy)(z)\;\mbox{ implies that }\;z=0,
\end{equation}
which amounts the kernel condition ${\rm ker}D^*F(\ox,\oy)=\{0\}$.

\section{The Generalized Newton Algorithm}
\setcounter{equation}{0}

This section presents the main contribution of the paper: a new
generalized Newton method for nonsmooth equations, which is based
on graphical derivatives. The section consists of three parts. In
Subsection~3.1 we precisely describe the algorithm and justify its
well-posedness/solvability. Subsection~3.2 contains a local superlinear
convergence result under appropriate assumptions. Finally,
in Subsection~3.3 we establish a global convergence result of the
Kantorovich type for our generalized Newton algorithm.

\subsection{Description and Justification of the Algorithm}

Keeping in mind the classical scheme of the smooth Newton method
in \eqref{eq:iteration}, \eqref{eq:linearized} and taking into
account the graphical derivative representation of
Proposition~\ref{prop:derivprop}(f), we propose an extension of
the Newton equation \eqref{eq:linearized} to nonsmooth mappings
given by:
\begin{equation}\label{eq:gennewteq}
-H(x^k)\in DH(x^k)(d^k),\quad k=0,1,2,\ldots.
\end{equation}
This leads us to the following generalized Newton algorithm to solve
\eqref{eq:ne}:\vspace*{0.05in}
\begin{Algorithm}\label{alg1} {\bf (generalized Newton's method).}

{\bf Step~0:} Choose a starting point $x^0\in\R^n$.\vspace*{0.05in}

{\bf Step~1:} Check a suitable termination criterion.\vspace*{0.05in}

{\bf Step~2:} Compute $d^k\in\R^n$ such that \eqref{eq:gennewteq}
holds.\vspace*{0.05in}

{\bf Step~3:} Set $x^{k+1}:=x^k+d^k$, $k\leftarrow k+1$, and go to Step~{\rm 1}.
\end{Algorithm}

\noindent The proposed Algorithm~\ref{alg1} does not require a
priori any assumptions on the underlying mapping $H\colon\R^n\to
\R^n$ in \eqref{eq:ne} besides its continuity, which is the
standing assumption in this paper. Other assumptions are imposed
below to justify the well-posedness and (local and global)
convergence of the algorithm. Observe that
Proposition~\ref{prop:derivprop}(c,d) ensures that
Algorithm~\ref{alg1} reduces to scheme \eqref{b1} in the
$B$-differentiable Newton method provided that $H$ is
directionally differentiable and locally Lipschitzian around the
solution point in question. In Section~5 we consider in detail
relationships with known results for the $B$-differentiable Newton
method, while Section~4 compares Algorithm~\ref{alg1} and the
assumptions made with the corresponding semismooth versions in the
framework of \eqref{c1}.\vspace*{0.05in}

To proceed further, we need to make sure that the generalized
Newton equation \eqref{eq:gennewteq} is {\em solvable}, which is a
major part of the well-posedness of Algorithm~\ref{alg1}. The next
proposition shows that an appropriate assumption to ensure the
solvability of \eqref{eq:gennewteq} is {\em metric regularity}.

\begin{Proposition}\label{lem:solv} {\bf (solvability of the generalized
Newton equation).}  Assume that $H\colon\R^n\to\R^n$ is metrically
regular around $\ox$ with $\oy=H(\ox)$ in \eqref{mr}, i.e., we
have $\ker D^*H(\ox)=\{0\}$. Then there is a constant $\ve>0$ such
that for all $x\in B_\ve(\ox)$ the equation
\begin{equation}\label{eq:gennewt}
   -H(x)\in DH(x)(d)
\end{equation}
admits a solution $d\in\R^n$. Furthermore, the set $S(x)$ of solutions to
\eqref{eq:gennewt} is computed by
\begin{equation}\label{s}
   S(x)=\Limsup_{t\dn 0,\,h\to-H(x)}\frac{H^{-1}\big(H(x)+th\big)-x}{t}\ne\emp.
\end{equation}
\end{Proposition}
{\bf Proof.} By the assumed metric regularity \eqref{mr} of $H$ we
find a number $\mu>0$ and neighborhoods $U$ of $\ox$ and $V$ of
$H(\ox)$ such that
\begin{equation*}
{\rm dist}\big(x;H^{-1}(y)\big)\le\mu\,{\rm dist}(y;H(x)\big)\;
\mbox{ for all }\;x\in U\;\mbox{ and }\;y\in V.
\end{equation*}
Pick now an arbitrary vector $x\in U$ and select sequences
$h_k\to-H(x)$ and $t_k\dn 0$ as $k\to\infty$. Suppose with no loss
of generality that $H(x)+t_k h_k\in V$ for all $k\in\N$. Then we
have
\begin{equation*}
{\rm dist}\big(x;H^{-1}(H(x)+t_k h_k)\big)\le\mu t_k\|h_k\|,\quad
k\in\N,
\end{equation*}
and hence there is a vector $u_k\in H^{-1}(H(x)+t_k h_k)$ such
that $\|u_k-x\|\le\mu t_k\|h_k\|$ for all $k\in\N$. This shows
that the sequence $\{\|u_k-x\|/t_k\}$ is bounded, and thus it
contains a subsequence that converges to some element $d\in\R^n$.
Passing to the limit as $k\to\infty$ and recalling the definitions
of the outer limit \eqref{pk} and of the tangent cone \eqref{tan},
we arrive at
\begin{equation*}
\big(d,-H(x)\big)\in\Limsup_{t\dn 0}\frac{\gph
H-\big(x,H(x)\big)}{t}= T\big((x,H(x));\gph H\big),
\end{equation*}
which justifies the desired inclusion \eqref{eq:gennewt}. The
solution representation \eqref{s} follows from \eqref{gd} and
Proposition~\ref{prop:derivprop}(b) in the case of single-valued
mappings, since
\begin{equation*}
S(x)=DH(x)^{-1}\big( \hspace{-1.5mm} -H(x) \hspace{-1.5mm} \big)
\end{equation*}
due to \eqref{eq:gennewt}. This completes the proof of the proposition. $\h$

\subsection{Local Convergence}

In this subsection we first formulate major assumptions of our
generalized Newton method and then show that they ensure the
superlinear local convergence of Algorithm~\ref{alg1}.
\begin{itemize}
\item[\bf (H1)] There exist a constant $C>0$, a neighborhood $U$ of $\bar x$,
and a neighborhood $V$ of the origin in $\R^n$ such that the following holds:\\
For all $x\in U$, $z\in V$, and for any $d\in\R^n$ with $-H(x)\in
DH(x)(d)$ there is a vector $w\in\tilde DH(x)(z)$ such that
\begin{equation*}
C\|d-z\|\le\|w+H(x)\|+o(\|x-\bar x\|).
\end{equation*}
\item[\bf (H2)] There exists a neighborhood $U$ of $\bar x$ such that for all
$v\in\tilde DH(x)(\bar x-x)$ we have
\begin{equation*}
\|H(x)-H(\bar x)+v\|=o(\|x-\bar x\|).
\end{equation*}
\end{itemize}
A detailed discussion of these two assumptions and sufficient
conditions for their fulfillment are given in Section~4. Note that
assumption (H2) means, in the terminology of
\cite[Definition~7.2.2]{FaP 03} focused on locally Lipschitzian
mappings $H$, that the family $\{\Tilde D H(x)\}$ provides a {\em
Newton approximation scheme} for $H$ at $\ox$.\vspace*{0.05in}

Now we establish our principal local convergence result that makes
use of the major assumptions (H1) and (H2) together with metric
regularity.

\begin{Theorem}\label{th:localconv} {\bf (superlinear local convergence of
the generalized Newton method).} Let $\ox\in\R^n$ be a solution to
\eqref{eq:ne} for which the underlying mapping
$H\colon\R^n\to\R^n$ is metrically regular around $\ox$ and
assumptions {\rm (H1)} and {\rm (H2)} are satisfied. Then there is
a number $\ve>0$ such that for all $x^0\in B_\ve(\ox)$ the
following assertions hold:

{\bf (i)} Algorithm~{\rm\ref{alg1}} is well defined and generates a
sequence $\{x^k\}$ converging to $\ox$.

{\bf (ii)} The rate of convergence $x^k\to\ox$ is at least superlinear.
\end{Theorem}
{\bf Proof.} To justify (i), pick $\ve>0$ such that assumptions
(H1) and (H2) hold with $U:=B_\ve(\ox)$ and $V:=B_\ve(0)$ and such
that Proposition~\ref{lem:solv} can be applied. Then we choose a
starting point $x^0\in B_\ve(\ox)$ and conclude by
Proposition~\ref{lem:solv} that the subproblem
\begin{equation*}
-H(x^0)\in DH(x^0)(d)
\end{equation*}
has a solution $d^0$. Thus the next iterate $x^1:=x^0+d^0$ is well
defined. Let further $z^0:=\ox-x^0$ and get $\|z^0\|\le\ve$ by the
choice of the starting point $x^0$. By assumption (H1), find a
vector $w^0\in\Tilde D H(x^0)(z^0)$ such that
\begin{equation*}
C\|x^1-\ox\|=C\|(x^1-x^0)-(\ox-x^0)\|=C\|d^0-z^0\|
\le\|w^0+H(x^0)\|+o(\|x^0-\ox\|).
\end{equation*}
Taking this into account and employing assumption (H2), we get the
relationships
\begin{eqnarray*}
C\|x^1-\ox\|&\le&\|w^0+H(x^0)\|+o(\|x^0-\ox\|)\\
&=&\|H(x^0)-H(\ox)+w^0\|+o(\|x^0-\ox\|)\\
&=&o(\|x^0-\ox\|)\\
&\le&\tfrac{C}{2}\|x^0-\ox\|,
\end{eqnarray*}
which imply that $\|x^1-\ox\|\le\frac{1}{2}\|x^0-\ox\|$. The
latter yields, in particular, that $x^1\in B_\ve(\ox)$. Now
standard induction arguments allow us to conclude that the
iterative sequence $\{x^k\}$ generated by Algorithm~\ref{alg1} is
well defined and converges to the solution $\ox$ of \eqref{eq:ne}
with at least a linear rate. This justifies assertion (i) of the
theorem.\vspace*{0.05in}

Next we prove assertion (ii) showing that the convergence
$x^k\to\ox$ is in fact {\em superlinear} under the validity of
assumption (H2). To proceed, we basically follow the proof of
assertion (i) and construct by induction sequences $\{d^k\}$
satisfying
\begin{equation*}
-H(x^k)\in DH(x^k)(d^k)\;\mbox{ for all }\;k\in\N,
\end{equation*}
$\{z^k\}$ with $z^k:=\ox-x^k$, and $\{w^k\}$ with $w^k\in\Tilde D
H(x^k)(z^k)$ such that
\begin{equation*}
C\|x^{k+1}-\ox\|\le\|w^k+H(x^k)\|+o(\|x^k-\ox\|),\quad k\in\N.
\end{equation*}
Applying then assumption (H2) gives us the relationships
\begin{equation*}
C\|x^{k+1}-\ox\|\le\|H(x^k)-H(\ox)+w^k\|+o(\|x^k-\ox\|)=o(\|x^k-\ox\|),
\end{equation*}
which ensure the superlinear convergence of the iterative sequence
$\{x^k\}$ to the solution $\ox$ of \eqref{eq:ne} and thus complete
the proof of the theorem. $\h$

\subsection{Global Convergence}

Besides the local convergence in the classical Newton method based
on suitable assumptions imposed at the (unknown) solution of the
underlying system of equations, there are global (or semi-local)
convergence results of the Kantorovich type \cite{KaA 64} for
smooth systems of equations which show that, under certain
conditions at the starting point $x^0$ and a number of assumptions
to hold in a suitable region around $x^0$, Newton's iterates are
well defined and converge to a solution belonging to this region;
see \cite{DeS 83,KaA 64} for more details and references. In the
case of nonsmooth equations \eqref{eq:ne} results of the
Kantorovich type were obtained in \cite{QiS 93,Rob 94} for the
corresponding versions of Newton's method. Global convergence
results of different types can be found in, e.g., \cite{DFK
96,FaP03,hpr,Pan 91} and their references. \vspace*{0.05in}

Here is a global convergence result for our generalized Newton
method to solve \eqref{eq:ne}.

\begin{Theorem}\label{gc} {\bf (global convergence of the
generalized Newton method).} Let $x^0$ be a starting
point of Algorithm~{\rm\ref{alg1}}, and let
\begin{equation}\label{om}
   \O:=\big\{x\in\R^n\big|\;\|x-x^0\|\le r\big\}
\end{equation}
with some $r>0$. Impose the following assumptions:
\begin{itemize}
\item[\bf(a)] The mapping $H\colon\R^n\to\R^n$ in \eqref{eq:ne} is
metrically regular on $\O$ with modulus $\mu>0$, i.e., it is
metrically regular around every point $x\in\O$ with the same
modulus $\mu$.

\item[\bf(b)] The set-valued map  $DH(x)(z)$ uniformly on $\O$
converges to $\{0\}$ as $z\to 0$ in the sense that: for all
$\ve>0$ there is $\delta>0$ such that
\begin{equation*}
   \|w\|\le\ve\;\mbox{ whenever }\;w\in DH(x)(z),\;\|z\|\le\dd,\;
   \mbox{ and }\;x\in\O.
\end{equation*}

\item[\bf(c)] There is $\al\in(0,1/\mu)$ such that
\begin{equation}\label{g1}
   \mu\|H(x^0)\|\le r(1-\al\mu)
\end{equation}
and for all $x,y\in\O$ we have the estimate
\begin{equation}\label{g2}
   \|H(x)-H(y)-v\|\le\al\|x-y\|\;\mbox{ whenever }\;v\in DH(x)(y-x).
\end{equation}
\end{itemize}
Then Algorithm~{\rm\ref{alg1}} is well defined, the sequence of
iterates $\{x^k\}$ remains in $\O$ and converges to a solution
$\ox\in\O$ of \eqref{eq:ne}. Moreover, we have the error estimate
\begin{equation}\label{er}
   \|x^k-\ox\|\le\frac{\al\mu}{1-\al\mu}\|x^k-x^{k-1}\|\;
   \mbox{ for all }\;k\in\N.
\end{equation}
\end{Theorem}

\noindent
{\bf Proof.} The metric regularity assumption (a) allows us to
employ Proposition~\ref{lem:solv} and, for any $x\in\O$ and
$d\in\R^n$ satisfying the inclusion $-H(x)\in DH(x)(d)$, to find
sequences of $h_k\to-H(x)$ and $t_k\dn 0$ as $k\to\infty$ such
that
\begin{equation*}
   \|d\|=\lim_{k\to\infty}\Big\|\frac{H^{-1}\big(H(x)+t_k
   h_k\big)-x}{t_k}\Big\|\le\lim_{k\to\infty}\mu\|h_k\|=\mu\|H(x)\|.
\end{equation*}
In view of assumption \eqref{g1} in (c) and the iteration
procedure of the algorithm, this implies
\begin{equation*}
   \|x^1-x^0\|=\|d^0\|\le\mu\|H(x^0)\|\le r(1-\al\mu),
\end{equation*}
which ensures that $x^1\in\O$ due the form of $\O$ in \eqref{om}
and the choice of $\al$. Proceeding further by induction, suppose
that $x^1,\ldots,x^k\in\O$ and get the relationships
\begin{align*}
   \|x^{k+1}-x^k\|&=\|d^k\|\le\mu\|H(x^k)\|\\
   &\le\mu\|H(x^k)-H(x^{k-1})+H(x^{k-1})\|\\
   &\le\al\mu\|x^k-x^{k-1}\|\quad\Big(\text{using \eqref{g2} and }
   -H(x^{k-1})\in DH(x^{k-1})(x^k-x^{k-1})\Big)\\
   &\le(\al\mu)^k\|x^1-x^0\|\le r(\al\mu)^k(1-\al\mu),
\end{align*}
which imply the estimates
\begin{equation*}
   \|x^{k+1}-x^0\|\leq\sum_{j=0}^k\|x^{j+1}-x^j\|\le\sum_{j=0}^k
   r(\al\mu)^j(1-\al\mu)\le r
\end{equation*}
and hence justify that $x^{k+1}\in\O$. Thus all the iterates
generated by Algorithm~\ref{alg1} remain in $\O$. Furthermore, for
any natural numbers $k$ and $m$, we have
\begin{equation*}
   \|x^{k+m+1}-x^k\|\le\sum_{j=k}^{k+m}\|x^{j+1}-x^j\|\le\sum_{j=k}^{k+m}
   r(\al\mu)^j(1-\al\mu)\le r(\al\mu)^k,
\end{equation*}
which shows that the generated sequence $\{x^k\}$ is a Cauchy
sequence. Hence it converges to some point $\ox$ that obviously
belongs to the underlying closed set \eqref{om}.

To show next that $\ox$ is a solution to the original equation
\eqref{eq:ne}, we pass to the limit as $k\to\infty$ in the
iterative inclusion
\begin{equation}\label{it}
-H(x^k)\in DH(x^k)(x^{k+1}-x^k),\quad k\in\N.
\end{equation}
It follows from assumption (b) that $\lim_{k\to\infty}H(x^k)=0$.
The continuity of $H$ then implies that $H(\ox)=0$, i.e., $\ox$ is
a solution to \eqref{eq:ne}.

It remains to justify the error estimate \eqref{er}. To this end,
first observe by \eqref{g1} that
\begin{equation*}
\|x^{k+m+1}-x^k\|\le\sum_{j=k}^{k+m}\|x^{j+1}-x^j\|\le\sum_{j=0}^{m}
(\al\mu)^{j+1}\|x^k-x^{k-1}\|\le\frac{\al\mu}{1-\al\mu}\|x^k-x^{k-1}\|
\end{equation*}
for all $k,m\in\N$. Passing now to the limit as $m\to\infty$, we
arrive at \eqref{er} thus completes the proof of the theorem. $\h$

\section{Discussion of Major Assumptions and
Comparison with Semismooth Newton Methods}
\setcounter{equation}{0}

In this section we pursue a twofold goal: to discuss the major
assumptions made in Section~3 and to compare our generalized
Newton method based on graphical derivatives with the semismooth
versions of the generalized Newton method developed in \cite{Qi
93,QiS 93}. As we will see from the discussions below, these two
aims are largely interrelated.

Let us begin with sufficient conditions for metric regularity in
terms of the constructions used in the semismooth versions of the
generalized Newton method. Given a locally Lipschitz continuous
vector-valued mapping $H\colon\R^n\to\R^m$, we have by the
classical Rademacher theorem that the set of points
\begin{equation}\label{R}
   S_H:=\{x\in\R^n\big|\;H\;\mbox{ is differentiable at }\;x\big\}
\end{equation}
is of full Lebesgue measure in $\R^n$. Thus for any mapping
$H\colon\R^n\to\R^m$ locally Lipschitzian around $\ox$ the set
\begin{equation}\label{eq:B-sub}
   \partial_B H(\ox):=\Big\{\lim_{k\to\infty}H'(x^k)\Big|\;\exists\,\{x^k\}
   \subset S_H\;\mbox{ with }\;x^k\to\ox\Big\}
\end{equation}
is nonempty and obviously compact in $\R^m$. It was introduced in
\cite{shor} for $m=1$ as the set of ``almost-gradients" and then
was called in \cite{Qi 93} the $B$-{\em subdifferential} of $H$ at
$\ox$. Clarke's {\em generalized Jacobian} \cite{Cla 83} of $H$ at
$\ox$ is defined by the convex hull
\begin{equation}\label{c}
   \partial_C H(\ox):=\mbox{co}\big\{\partial_B H(\ox)\big\}.
\end{equation}
We also make use of the {\em Thibault derivative/limit set}
\cite{Thi 82} (called sometimes the ``strict graphical derivative"
\cite{rw}) of $H$ at $\ox$ defined by
\begin{equation}\label{t}
D_TH(\ox)(z):= \Limsup_{x\to\ox\atop t \downarrow 0}\frac{H(x+t
z)-H(x)}{t},\quad z\in\R^n.
\end{equation}
Observe the known relationships \cite{kk,Thi 82} between the above
derivative sets
\begin{equation}\label{Eq:RelThibault}
\partial_B H(\ox)z\subset D_T H(\ox)(z)\subset\partial_C
H(\ox)z,\quad z\in\R^n.
\end{equation}

The next result gives a sufficient condition for metric regularity
of Lipschitzian mappings in terms of the Thibault derivative
\eqref{t}. It can be derived from the coderivative
characterization of metric regularity
\eqref{eq:coderivativecriterion}, while we give here a direct
independent proof.

\begin{Proposition}\label{lem:nonsing} {\bf (sufficient condition for
metric regularity in terms of Thibault's derivative).} Let
$H\colon\R^n\to\R^n$ be locally Lipschitzian around $\ox$, and let
\begin{equation}\label{smr}
0\notin D_T H(\ox)(z)\;\mbox{ whenever }\;z\ne 0.
\end{equation}
Then the mapping $H$ is metrically regular around $\ox$.
\end{Proposition}

\noindent
{\bf Proof.} Kummer's inverse function theorem
\cite[Theorem~1.1]{k} ensures that condition \eqref{smr} implies
(actually is equivalent to) the fact that there are neighborhoods
$U$ of $\ox$ and $V$ of $H(\ox)$ such that the mapping $H\colon
U\to V$ is one-to-one with a locally Lipschitzian inverse
$H^{-1}\colon V\to U$. Let $\mu>0$ be a Lipschitz constant of
$H^{-1}$ on $V$. Then for all $x\in U$ and $y\in V$ we have the
relationships
\begin{eqnarray*}
   \text{dist}\big(x;H^{-1}(y)\big)&=&\|x-H^{-1}(y)\|\\
   &=&\|H^{-1}\big(H(x)\big)-H^{-1}(y)\|\\
   &\le&\mu\|H(x)-y\|\\
   &=&\mu\,\text{dist}\big(y;H(x)\big),
\end{eqnarray*}
which thus justify the metric regularity of $H$ around $\ox$.
$\h$\vspace*{0.05in}

To proceed further with sufficient conditions for the validity of
our assumption (H1), we first introduce the notion of directional
boundedness.

\begin{Definition}\label{dir-bound} {\bf (directional
boundedness).} A mapping $H\colon\R^n\to\R^m$ is said to be {\sc
directionally bounded} around $\ox$ if
\begin{equation}\label{Eq:DirBounded}
\limsup_{t\downarrow
0}\left\|\frac{H(x+tz)-H(x)}{t}\right\|<\infty
\end{equation}
for all $x$ near $\ox$ and for all $z\in\R^n$.
\end{Definition}

It is easy to see that if $H$ is either directionally
differentiable around $\ox$ or locally Lip\-schitzian around this
point, then it is directionally bounded around $\ox$. The
following example shows that the converse does not hold
in general.

\begin{Example}\label{dir1} {\bf (directional bounded mappings
may not be directionally differentiable).}
{\rm Define a real function $H\colon\R\to\R$ by
\begin{equation*}
   H(x):=\left\{\begin{array}{ll}
   x\sin\big(\frac{1}{x}\big)&\text{if }\;x \ne 0,\\
   0&\text{if }\;x=0.
\end{array}\right.
\end{equation*}
It is easy to see that this function is not directionally
differentiable at $\bar x=0$. However, it is directionally bounded
around $\bar x$. Indeed, for any $x\ne 0$ near $\bar x$ condition
\eqref{Eq:DirBounded} holds because $H$ is simply differentiable
at $x\ne 0$. For $x=0$ we have
\begin{equation*}
   \limsup_{t\downarrow 0}\Big|\frac{H(tz)-H(0)}{t}\Big|=\limsup_{
   t\downarrow 0}\frac{|H(tz)|}{t}=\limsup_{t\downarrow 0}\Big|z\sin
   \Big(\frac{1}{tz}\Big)\Big|=|z|<\infty.
\end{equation*}}
\end{Example}

The next proposition and its corollary present verifiable
sufficient conditions for the fulfillment of assumption (H1).

\begin{Proposition}\label{prop:suffH1} {\bf (assumption (H1) from
metric regularity).} Let $H\colon\R^n\to\R^n$, and let $\ox$ be a
solution to \eqref{eq:ne}. Suppose that $H$ is metrically regular
around $\ox$ $($i.e., $\ker D^*H(\ox)=0)$, that it is
directionally bounded and one-to-one around this point. Then
assumption {\rm (H1)} is satisfied.
\end{Proposition}

\noindent
{\bf Proof.} Recall that the metric regularity of $H$ around $\ox$
is equivalent to the condition $\ker D^*H(\ox)=\{0\}$ by the
coderivative criterion \eqref{eq:coderivativecriterion}. Let
$U\subset\R^n$ be a neighborhood of $\ox$ such that $H$ is
metrically regular and one-to-one on $U$. Choose further a
neighborhood $V\subset\R^n$ of $H(\ox)=0$ from the definition of
metric regularity of $H$ around $\ox$. Then pick $x\in U$, $z\in
V$ and an arbitrary direction $d\in\R^n$ satisfying $-H(x)\in
DH(x)(d)$. Employing now Proposition~\ref{lem:solv}, we get
\begin{equation*}
   d\in\Limsup_{h\to-H(x),\;t\downarrow
   0}\frac{H^{-1}\big(H(x)+th\big)-x}{t}.
\end{equation*}
By  the local single-valuedness of $H^{-1}$ and the metric
regularity of $H$ around $\ox$ there exists a number $\mu>0$ such
that
\begin{equation*}
   \left\|\frac{H^{-1}(H(x)+th)-x}{t}-z\right\|\le
   \mu\left\|\frac{H(x)+th-H(x+tz)}{t}\right\|=\mu\left\|\frac{H(x+tz)-H(x)}
   {t}-h\right\|
\end{equation*}
for all $t>0$ sufficiently small. It follows that
\begin{equation*}
   \|d-z\|\le\limsup_{t\downarrow 0\atop h\to-H(x)}\left\|
   \frac{H^{-1}\big(H(x)+th\big)-x}{t}-z\right\|\le\mu\limsup_{t\downarrow
   0\atop h\to-H(x)}\left\|\frac{H(x+tz)-H(x)}{t}-h\right\|<\infty
\end{equation*}
by the directional boundedness of $H$ around $\ox$. The
boundedness of the family
\begin{equation*}
   \Big\{v(t):=\frac{H(x+tz)-H(x)}{t}\Big\},\quad t \downarrow 0,
\end{equation*}
allows us to select a sequence $t_k\downarrow 0$ such that $v(t_k)\to
w$ for some $w\in\R^n$. By passing to the limit above as
$k\to\infty$ and employing Definition~\ref{def:restrder} we get
that
\begin{equation*}
   w\in\Tilde DH(x)(z)\quad{\rm and}\quad
   \frac{1}{\mu}\|d-z\|\le\|w+H(x)\|,
\end{equation*}
which completes the proof of the proposition. $\h$\vspace*{0.05in}

\begin{Corollary}\label{suffH1} {\bf (sufficient
conditions for (H1) via Thibault's derivative).} Let $\ox$ be a
solution to \eqref{eq:ne}, where $H\colon\R^n\to\R^n$ is locally
Lipschitzian around $\ox$ and such that condition \eqref{smr}
holds, which is automatic when ${\rm det}\,A\ne 0$ for all
$A\in\partial_C H(\ox)$. Then {\rm(H1)} is satisfied with $H$
being both metrically regular and one-to-one around $\ox$.
\end{Corollary}

\noindent
{\bf Proof.} Indeed, both metric regularity and bijectivity of $H$
around $\ox$ assumed in Proposition~\ref{prop:suffH1} follow from
Proposition~\ref{lem:nonsing} and its proof. Nonsingularity of all
$A\in\partial_C H(\ox)$ clearly implies \eqref{smr} by the second
inclusion in \eqref{Eq:RelThibault}. $\h$\vspace*{0.05in}

Note that other conditions ensuring the fulfillment of assumption
(H1) for Lipschitzian and non-Lipschitzian mappings
$H\colon\R^n\to\R^n$ can be formulated in terms of Warga's {\em
derivate containers} by \cite[Theorems~1 and 2]{w} on ``fat
homeomorphisms" that also imply the metric regularity and
one-to-one properties of $H$.\vspace*{0.05in}

Next we proceed with the discussion of assumption (H2) and
present, in particular, sufficient conditions for their
fulfillment via semismoothness. First observe the following.

\begin{Proposition}\label{prop:graphgen} {\bf (relationship
between graphical derivative and generalized Jacobian).} Let
$H\colon\R^n\to\R^m$ be locally Lipschitzian around $\ox$. Then we
have
\begin{equation}\label{c2}
    DH(\ox)(z)\subset\partial_C H(\ox)z\;\mbox{ for all }\;z\in\R^n.
\end{equation}
\end{Proposition}

\noindent
{\bf Proof.} Pick $w\in DH(\ox)(z)$ and get by
Proposition~\ref{prop:derivprop}(c) and
Definition~\ref{def:restrder} a sequence of $t_k\dn 0$ as
$k\to\infty$ such that
\begin{equation}\label{w}
   w=\lim_{k\to\infty}\frac{H(\ox+t_kz)-H(\ox)}{t_k}.
\end{equation}
It follows from \cite[Proposition~2.6.5]{Cla 83} that
\begin{equation*}
   \frac{H(\ox+t_kz)-H(\ox)}{t_k}\in{\rm co}\big\{\partial_C
   H[\ox,\ox+t_k z]\big\}z\;\mbox{ for all }\;k\in\N.
\end{equation*}
Applying to the latter the classical Carath\'eodory theorem, we
find scalars $\gg^k_i\in[0,t_k]$, $\lm^k_i\in[0,1]$ and matrices
$A^k_i\in\partial_C H(\ox+\gg^k_i z)$ for $i=1,\ldots,m+1$ such
that
\begin{equation*}
   \frac{H(\ox+t_kz)-H(\ox)}{t_k}=\Big[\sum_{i=1}^{m+1}\lambda_i^kA_i^k\Big]
   z\quad\text{and}\quad\sum_{i=1}^{m+1}\lambda_i^k=1\;\mbox{ for all }\;
   k\in\N.
\end{equation*}
Due to the boundedness of the sequences $\{\lm^k_i\}_{k\in\N}$,
the convergence $\ox+\gg^k_i z\to\ox$ as $k\to\infty$ for all
$i=1,\ldots,m+1$, and the outer/upper semicontinuity of the
mapping $x\mapsto\partial_C H(x)$ proved in
\cite[Proposition~2.6.2]{Cla 83} we have that the sequences
$\{A^k_i\}$ are bounded as well. Hence there are subsequences of
these sequences (without relabelling), scalars $\lm_i\in[0,1]$,
and matrices $A_i$ as $i=1,\ldots,m+1$ such that
\begin{equation*}
   \lm^k_i\to\lm_i,\quad\sum_{i=1}^{m+1}\lm_i=1,\;\mbox{ and
   }\;A^k_i\to A_i\in\partial_C H(\ox)\;\mbox{ as }\;k\to\infty.
\end{equation*}
By \eqref{w} and the subsequent relationships therein, we get
\begin{equation*}
   w=\lim_{k\to\infty}\Big[\sum_{i=1}^{m+1}\lambda_i^kA_i^k\Big]z =
   \Big[\sum_{i=1}^{m+1}\lambda_iA_i\Big]z\in{\rm
   co}\big\{\partial_{C}H(\ox)\big\}z=\partial_{C}H(\ox)z
\end{equation*}
and thus complete the proof of the proposition.
$\h$\vspace*{0.05in}

Inclusion \eqref{c2}---which may be strict as illustrated by
Example~\ref{ex:counterex2} below---shows that our generalized
Newton Algorithm~\ref{alg1} based on the graphical derivative
provides in the case of Lipschitz equations \eqref{eq:ne} a more
accurate choice of the iterative direction $d^k$ via
\eqref{eq:gennewteq} in comparison with the iterative relationship
\begin{equation}\label{semi}
   -H(x^k)\in\partial_C H(x^k)d^k,\quad k=0,1,2,\ldots,
\end{equation}
used in the semismooth Newton method \cite{QiS 93} and related
developments \cite{kk, k88} based on the generalized Jacobian. If
in addition to the assumptions of Proposition~\ref{prop:graphgen}
the mapping $H$ is directionally differentiable at $\ox$, then
$DH(\ox)(z)=\{H'(\ox;z)\}$ by
Proposition~\ref{prop:derivprop}(c,d). Thus in this case we have
from Proposition~\ref{prop:graphgen} that for any $z\in\R^n$ there
is $A\in\partial_C H(\ox)$ such that $H'(\ox;z)=Az$, which
recovers a well-known result from \cite[Lemma~2.2]{QiS
93}.\vspace*{0.05in}

The following example shows that the converse inclusion in
Proposition~\ref{prop:graphgen} is not satisfied in general even
with the replacement of the set $DH(\ox)(z)$ in \eqref{c2} by its
convex hull co$DH(\ox)(z)$ in the case of real functions.
Furthermore, the same holds if we replace the generalized Jacobian
in \eqref{c2} by the smaller $B$-subdifferential $\partial_B
H(\ox)$ from \eqref{eq:B-sub}.

\begin{Example}\label{ex:counterex2} {\bf (graphical derivative is
strictly smaller than $B$-subdifferential and generalized
Jacobian).} {\rm Consider the simplest nonsmooth convex function
$H(x)=|x|$ on $\R$. In this case $\partial_B H(0)=\{-1,1\}$ and
$\partial_C H(0)=[-1,1]$. Thus
\begin{equation*}
   \partial_B H(0)z=\{-1,1\}\;\mbox{ and }\;\partial_C
   H(0)z=[-1,1]\;\mbox{ for }\;z=1.
\end{equation*}
Since $H(x)=|x|$ is locally Lipschitzian and directionally
differentiable, we have
\begin{equation*}
   DH(0)(z)=\big\{H'(0;z)\big\}=|z|=\{1\}\;\mbox{ for }\;z=1.
\end{equation*}
Hence it gives the relationships
\begin{equation*}
   DH(0)(z)={\rm co}\big\{DH(0)(z)\big\}\subset\partial_B
   H(0)z\subset\partial_C H(0)z,
\end{equation*}
where both inclusions are strict. Observe also the difference
between the convexification of the graphical derivative and of the
coderivative; in the latter case we have equality \eqref{C2}.}
\end{Example}

As mentioned in Section~1, there is an improvement \cite{Qi 93} of
the iterative procedure \eqref{semi} with the replacement the
generalized Jacobian therein by the $B$-subdifferential
\begin{equation}\label{semi1}
   -H(x^k)\in\partial_B H(x^k)d^k,\quad k=0,1,2,\ldots.
\end{equation}
Note that, along with obvious advantages of version \eqref{semi1}
over the one in \eqref{semi}, in some settings it is easier to
deal with the generalized Jacobian than with its
$B$-subdifferential counterpart due to much better calculus and
convenient representations for $\partial_C H(\ox)$ in comparison
with the case of $\partial_B H(\ox)$, which does not even reduce
to the classical subdifferential of convex analysis for simple
convex functions as, e.g., $H(x)=|x|$. A remarkable common feature
for both versions in \eqref{semi} and \eqref{semi1} is the
efficient semismoothness assumption imposed on the underlying
mapping $H$ to ensure its local superlinear convergence. This
assumption, which unifies and labels versions \eqref{semi} and
\eqref{semi1} as the ``semismooth Newton method", is replaced in
our generalized Newton method by assumption (H2). Let us now
recall the notion of semismoothness and compare it with
(H2).\vspace*{0.05in}

A mapping $H\colon\R^n\to\R^m$, locally Lipschitzian and
directionally differentiable around $\ox$, is {\em semismooth} at
this point if the limit
\begin{equation}\label{eq:semismooth}
\lim_{h\to z,\;t\downarrow 0\atop A\in \partial_{C}H(\ox+th)}
\big\{Ah\big\}
\end{equation}
exists for all $z\in\R^n$; see \cite[Definition~7.4.2]{FaP 03}.
This notion was introduced in \cite{Mif 77} for real-valued
functions and then extended in \cite{QiS 93} to the vector
mappings for the purpose of applications to a nonsmooth Newton's
method. It is not hard to check \cite[Proposition~2.1]{QiS 93}
that the existence of the limit in \eqref{eq:semismooth} implies
the directional differentiability of $H$ at $\ox$ (but may not
around this point) with
\begin{equation*}
H'(\ox;z)=\lim_{h\to z,\;t\downarrow 0\atop A\in
\partial_C H(\ox+th)}\big\{Ah\big\}\;\mbox{ for all }\;z\in\R^n.
\end{equation*}
One of the most useful properties of semismooth mappings is the
following representation for them obtained in
\cite[Proposition~1]{PaQ 93}:
\begin{equation}\label{eq:semismooth1}
   \|H(\ox+z)-H(\ox)-Az\|=o(\|z\|)\;\mbox{ for all }\;z\to 0\;\mbox{
   and }\;A\in\partial_C H(\ox+z),
\end{equation}
which we exploit now to relate semismoothness to our assumption
(H2).

\begin{Proposition}\label{prop:suffH2} {\bf (semismoothness
implies assumption (H2)).} Let $H\colon\R^n\to\R^m$ be semismooth
at $\bar x$. Then assumption {\rm (H2)} is satisfied.
\end{Proposition}

\noindent
{\bf Proof.} Since any semismooth mapping is Lipschitz continuous
on a neighborhood $U$ of $\ox$, we have by
Proposition~\ref{prop:derivprop}(c) that
\begin{equation*}
\Tilde DH(x)(\ox-x)=DH(x)(\ox-x)\;\mbox{ for all }\;x\in U.
\end{equation*}
Proposition~\ref{prop:graphgen} yields therefore that
\begin{equation*}
\Tilde DH(x)(\ox-x)\subset\partial_C H(x)(\ox-x)\;\mbox{ whenever
}\;x\in U.
\end{equation*}
Given any $v\in\Tilde D H(x)(\ox-x)$ and using the latter
inclusion, find a matrix $A\in\partial_C H(x)$ such that
$v=A(\ox-x)$. Applying finally property \eqref{eq:semismooth1} of
semismooth mappings, we get
\begin{equation*}
\|H(x)-H(\ox)+v\|=\|H(x)-H(\ox)-A(x-\ox)\|=o(\|x-\ox \|)\;\mbox{
for all }\;x\in U,
\end{equation*}
which thus verifies (H2) and completes the proof of the
proposition. $\h$\vspace*{0.05in}

Note that the previous proposition actually shows that condition
\eqref{eq:semismooth1} implies (H2). The next result states that
the converse is also true, i.e., we have that assumption (H2) is
completely equivalent to \eqref{eq:semismooth1} for locally
Lipschitzian mappings.

\begin{Proposition}\label{prop:H2equiv} {\bf (equivalent description
of (H2)).} Let $H\colon\R^n\to\R^m$ be locally Lipschitzian around
$\ox$, and let assumption {\rm (H2)} hold with some neighborhood
$U$ therein. Then
\begin{equation}\label{semiB}
   \|H(x)-H(\ox)-A(x-\ox)\|=o(\|\ox-x\|)\;\mbox{ for all }\;x\in
   U\;\mbox{ and }\;A\in\partial_B H(x).
\end{equation}
Therefore assumption {\rm (H2)} is equivalent to
\eqref{eq:semismooth1}.
\end{Proposition}

\noindent
{\bf Proof.} Arguing by contradiction, suppose that \eqref{semiB}
is violated and find sequences $x_k\to\ox$, $A_k\in\partial_B
H(x_k)$ and a constant $\gg>0$ such that
\begin{equation*}
   \|H(x_k)-H(\ox)-A_k(x_k-\ox)\|\ge\gg\|\ox-x_k\|,\quad k\in\N.
\end{equation*}
By the Lipschitz property of $H$ and by construction
(\ref{eq:B-sub}) of the $B$-subdifferential there are points of
differentiability $u_k\in S_H$ close to $x_k$ with $H'(u_k)$
sufficiently close to $A_k$ satisfying
\begin{equation*}
   \|H(u_k)-H(\ox)-H'(u_k)(u_k-\ox)\|\ge\tfrac{\gg}{2}\|\ox-u_k\|,\quad
   k\in\N.
\end{equation*}
Then Proposition~\ref{prop:derivprop}(c,f) gives us the
representations
\begin{equation*}
   \Tilde DH(u_k)(\ox-u_k)=DH(u_k)(\ox-u_k)=-H'(u_k)(u_k-\ox)
\end{equation*}
for all $k\in\N$, which imply therefore that
\begin{equation*}
   \|H(u_k)-H(\ox)+v\|\ge\tfrac{\gg}{2}\|\ox-u_k\|\;\mbox{ whenever
   }\;v\in\Tilde DH(u_k)(\ox-u_k),\quad k\in\N.
\end{equation*}
This clearly contradicts assumption (H2) for $k$ sufficiently
large and thus ensures property \eqref{semiB}. The equivalence
between (H2) and \eqref{eq:semismooth1} follows now from the
implication (H2)$\Longrightarrow$\eqref{semiB} and the proof of
Proposition~\ref{prop:suffH2}. $\h$\vspace*{0.05in}

It is well known that, for the class of locally Lipschitzian and
directionally differentiable mappings, condition
\eqref{eq:semismooth1} is equivalent to the original definition of
semismoothness; see, e.g., \cite[Theorem~7.4.3]{FaP 03}.
Proposition~\ref{prop:H2equiv} above establishes the equivalence
of \eqref{eq:semismooth1} to our major assumption (H2) provided that
$H$ is locally Lipschitzian around the reference point while it
may {\em not} be directionally differentiable therein. In fact, it
follows from Example~\ref{exH1} that assumption (H2) may hold for
locally Lipschitzian functions, which are not directionally
differentiable and hence not semismooth. Let us now illustrate
that (H2) may also be satisfied for non-Lipschitzian mappings, in
which case it is {\em not} equivalent to property
\eqref{eq:semismooth1}.

\begin{Example}\label{non-lip} {\bf (assumption (H2) holds for
non-Lipschitzian one-to-one mappings).}  {\rm Consider the mapping
$H\colon\R^2\to\R^2$ defined by
\begin{equation}\label{non-lip1}
   H(x_1,x_2):=\Big(x_2\sqrt{|x_1|+|x_2|^3},x_1\Big)\;\mbox{ for
   }\;x_1, x_2 \in \R.
\end{equation}
It is easy to check that this mapping is
one-to-one around $(0,0)$. Focusing for definiteness on the
nonnegative branch of the mapping $H$, observe
that at any point $(x_1,x_2)\in\R^2$ with either $x_1, x_2>0$,
the classical Jacobian $J H(x_1,x_2)$ is computed by
\begin{eqnarray*}
JH(x_1,x_2)=\left[\begin{array}{c}\begin{array}{ll}\disp\frac{x_2}{2
   \sqrt{x_1+x^3_2}}\qquad\sqrt{x_1+x^3_2}+\disp\frac{3
   x^3_2}{2\sqrt{x_1+x^3_2}}\\\qquad1\qquad\qquad\qquad\qquad
   0\end{array}
   \end{array}\right].
\end{eqnarray*}
Setting $x_1=x_2^3$, we see that the first component
\begin{equation*}
   \frac{x_2}{2\sqrt{x_1+x_2^3}}=\frac{x_2}{2\sqrt{x_2^3+x_2^3}}
\end{equation*}
is unbounded when $x_1,x_2\dn 0$. This implies that the Jacobian
$JH(x_1,x_2)$ is unbounded around $(\ox_1,\ox_2)=(0,0)$, and hence
$H$ is not locally Lipschitzian around the origin.

Let us finally verify that the underlying assumption (H2) is
satisfied for the mapping $H$ in \eqref{non-lip1}. First assume that
$ x_1, x_2 > 0 $. Then we need to check that
$$
\begin{aligned}
   &\|H(x_1,x_2)-H(\ox_1,\ox_2)+JH(x_1,x_2)(-x_1,-x_2)\|\\
   &=\left|x_2\sqrt{x_1+x_2^3}-\frac{x_1x_2}{2\sqrt{x_1+x_2^3}}-x_2
   \sqrt{x_1+x_2^3}-\frac{3x_2^4}{2\sqrt{x_1+x_2^3}}\right|\\
   &=\left|\frac{x_1x_2}{2\sqrt{x_1+x_2^3}}+\frac{3x_2^4}{2\sqrt{x_1+x_2^3}}
   \right|=o\big(\sqrt{x_1^2+x_2^2}\big).
\end{aligned}
$$
The latter surely holds as $(x_1,x_2)\to (0,0)$ due to the
estimates
\begin{equation*}
   \frac{x_1x_2}{2\sqrt{x_1+x_2^3}\sqrt{x_1^2+x_2^2}}\le\frac{x_1}
   {\sqrt{x_1+x_2^3}}\le\sqrt{x_1},
\end{equation*}
\begin{equation*}
   \frac{3x_2^4}{2\sqrt{x_1+x_2^3}\sqrt{x_1^2+x_2^2}}\le\frac{3x_2^3}
   {2\sqrt{x_1+x_2^3}}\le 3x_2,
\end{equation*}
which thus justify the fulfillment of assumption (H2) in this
case. The other cases where $ x_1 > 0, x_2 \le 0 $ or $ x_1 < 0,
x_2 > 0 $ or $ x_1 < 0, x_2 \le 0 $ or, finally, $ x_1 = 0, x_2 $
arbitrary (here $H$ is not differentiable) can be treated in a
similar way.}
\end{Example}

To complete our discussion on the major assumptions in this
section, let us present an example of a locally Lipschitzian
function, which satisfies assumptions (H1) and (H2) being locally
one-to-one and metrically regular around the point in question
while not being directionally differentiable and hence not
semismooth at this point.

\begin{figure}[ht]
\begin{center}
\scalebox{0.6}{\input{Hfun.eepic}}
\end{center}
\caption{Construction of the mapping from Example \ref{exH1}:
Illustration}\label{Fig:Hfun}
\end{figure}
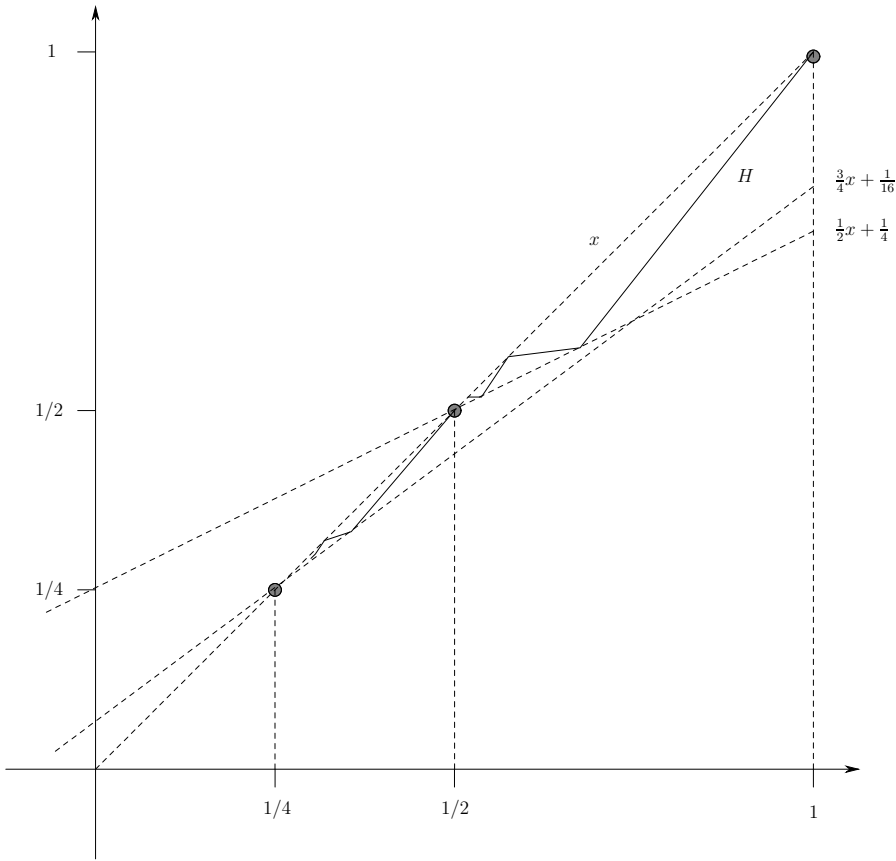

\begin{Example}\label{exH1} {\bf (non-semismooth but metrically
regular, Lipschitzian, and one-to-one functions satisfying (H1)
and (H2)).} {\rm We construct a function $H\colon[-1,1]\to\R$ in
the following way. First set $H(\ox):=0$ at $\ox=0$. Then define
$H$ on the interval $(1/2,1]$ staying between two lines
\begin{equation*}
\Big(1-\frac{1}{2}\Big)x+\frac{1}{4}\le H(x)\le x
\end{equation*}
in the following way: start from $(1,1)$ and let $H$ be continuous
piecewise linear when $x$ goes from 1 to 1/2 with the slope 1+1/4
and then with the slope $1/2-1/4$ alternatively until $x$ reaches
1/2. Consider further each interval $(2^{-k},2^{-(k-1)}]$ for
$k=2,3,\ldots$ and, starting from the point
$\big(2^{-(k-1)},2^{-(k-1)}\big)$, define $H$ to be continuous
piecewise linear with the corresponding slopes of either
$1+2^{-2k}$ or $1-2^{-k}-2^{-2k}$ staying between the two lines
\begin{equation}\label{sl}
\Big(1-\frac{1}{2^k}\Big)x+\frac{1}{2^{2k}}\le H(x)\le x.
\end{equation}
Thus we have constructed $H$ on the whole interval $[0,1]$; see
Figure~\ref{Fig:Hfun} for illustration. On the interval
$[-1,0]$, define the function $H$ symmetrically with respect to
the origin. Then it is easy to see that $H$ in continuous on
$[-1,1]$ and satisfies the following properties:
\begin{itemize}

\item $H$ is clearly Lipschitz continuous around $\ox=0$.

\item Since $H$ is continuous and monotone with a positive uniform
slope, it is one-to-one and metrically regular around $\ox$, which
directly follows, e.g., from the coderivative criterion
\eqref{eq:coderivativecriterion}. This ensures the fulfillment of
assumption (H1) by Proposition~\ref{prop:suffH1}.

\item To verify assumption (H2), fix $k\in\N$ and
$x\in(2^{-k},2^{-(k-1)}]$ and then pick any
\begin{equation*}
v\in
DH(x)(\ox-x)\subset\Big[1-\frac{1}{2^{k}}-\frac{1}{2^{2k}},1+\frac{1}{2^{2k}}\Big](\ox-x).
\end{equation*}
Since $\ox=0$, the latter implies that
\begin{equation*}
-\Big(1+\frac{1}{2^{2k}}\Big)x\le
v\le\Big(1-\frac{1}{2^k}-\frac{1}{2^{2k}}\Big)x
\end{equation*}
Thus we have by \eqref{sl} and simple computations that
\begin{equation*}
|H(x)-H(\ox)+v|\le\frac{1}{2^k}|x|+\frac{1}{2^{2k}}+\frac{1}{2^{2k}}=o\Big(\frac{1}{2^{k}}\Big)=
o(|x-\ox|),
\end{equation*}
which shows that assumption (H2) is satisfied. In  fact, it
follows from above that the latter value is
$O(2^{-2k})=O(\|x-\ox\|^2)$.

\item Let us finally check that $H$ is not directionally
differentiable at $x_k=2^{-k}$ for any $k\in\N$; therefore it is
not directionally differentiable around the reference point
$\ox=0$ and hence not semismooth at $\ox$. Indeed, this follows
directly from computing the graphical derivative by
\begin{equation*}
DH(x_k)(1)=\Big[1-\frac{1}{2^{k}},1\Big],\quad k\in\N,
\end{equation*}
which is not single-valued at $x_k$, and thus $H$ is not
directionally differentiable at $x_k$ due to
Proposition~\ref{prop:derivprop}(c,d).
\end{itemize}}
\end{Example}

\section{Application to the $B$-differentiable Newton  Method}
\setcounter{equation}{0}

In this section we present applications of the graphical
derivate-based generalized Newton method developed above to the
$B$-differentiable Newton method for nonsmooth equations
\eqref{eq:ne} originated by Pang \cite{Pan 90}.

Throughout this section, suppose that $H\colon\R^n\to\R^n$ is
locally Lipschitzian and directionally differentiable around the
reference solution $\ox$ to \eqref{eq:ne}.
Proposition~\ref{prop:derivprop}(c,d) yields in this setting that
the generalized Newton equation \eqref{eq:gennewteq} in our
Algorithm~\ref{alg1} reduces to
\begin{equation}\label{eq:B-diffnewt}
-H(x^k)=H'(x^k;d^k)
\end{equation}
with respect to the new search direction $d^k$ and that the new
iterate $x^{k+1}$ is computed by
\begin{equation}\label{eq:iter}
x^{k+1}:=x^k+d^k,\quad k=0,1,2,\ldots.
\end{equation}

Note that Pang's $B$-differentiable Newton method and its further
developments (see, e.g., \cite{FaP 03,hpr,Pan 91,Qi 93,QiS 93})
are based on Robinson's notion of the $B$(ouligand)-derivative
\cite{Rob 87} for nonsmooth mappings; hence the name. As was then
shown in \cite{sh}, the $B$-derivative of a locally Lipschitzian
mapping agrees with the classical directional derivative. Thus the
iteration scheme in Pang's $B$-differentiable method reduces to
\eqref{eq:B-diffnewt} and \eqref{eq:iter} in the Lipschitzian and
directionally differentiable case, and so we keep the original
name of \cite{Pan 90}.\vspace*{0.05in}

The next theorem shows what we get from applying our local
convergence result from Theorem~\ref{th:localconv} and the
subsequent analysis developed in Sections~3 and 4 to the
$B$-differentiable Newton method. This theorem employs an
equivalent description of assumption (H2) held in the setting
under consideration and the coderivative criterion
\eqref{eq:coderivativecriterion} for metric regularity of the
underlying Lipschitzian mapping $H$ ensuring the validity of
assumption (H1).

\begin{Theorem}\label{thm:Bdiffconv} {\bf (solvability and local
convergence of the $B$-differentiable Newton method via metric
regularity).} Let $H\colon\R^n\to\R^n$ be semismooth, one-to-one,
and metrically regular around a reference solution $\ox$ to
\eqref{eq:ne}, i.e.,
\begin{equation}\label{cod-cr}
0\in\partial\la z,H\ra(\ox)\Longrightarrow z=0.
\end{equation}
Then the $B$-differentiable Newton method \eqref{eq:B-diffnewt},
\eqref{eq:iter} is well defined $($meaning that equation
\eqref{eq:B-diffnewt} is solvable for $d^k$ as $k\in\N)$ and
converges at least superlinearly to the solution $\ox$.
\end{Theorem}

\noindent
{\bf Proof.} Since $H$ is locally Lipschitzian around $\ox$, the
coderivative criterion \eqref{eq:coderivativecriterion} is
equivalently written in form \eqref{cod-cr} via the limiting
subdifferential \eqref{sub} due to the scalarization formula
\eqref{Eq:scalarization}. Applying Theorem~\ref{th:localconv} to
the $B$-differentiable Newton method, we need to check that
assumptions (H1) and (H2) are satisfied in the setting under
consideration. Indeed, it follows from
Proposition~\ref{prop:H2equiv} and the discussion right after it
that (H2) is equivalent to the semismoothness for locally
Lipschitzian and directionally differentiable mappings. The
fulfillment of assumption (H1) is guaranteed by
Proposition~\ref{prop:suffH1}. $\h$\vspace*{0.05in}

More specific sufficient conditions for the well-posedness and
superlinear convergence of the $B$-differentiable Newton method
are formulated via of the Thibault derivative \eqref{t}.

\begin{Corollary}\label{Bsub-th} {\bf ($B$-differentiable Newton
method via Thibault's derivative).} Let $H\colon\R^n\to\R^n$ be
semismooth at the reference solution point $\ox$ of equation
\eqref{eq:ne}, and let condition \eqref{smr} be satisfied. Then
the $B$-subdifferential Newton method \eqref{eq:B-diffnewt},
\eqref{eq:iter} is well defined and converges superlinearly to the
solution $\ox$.
\end{Corollary}

\noindent
{\bf Proof.} Follows from Theorem~\ref{thm:Bdiffconv} and
Proposition~\ref{suffH1}. $\h$\vspace*{0.05in}

Observe by the second inclusion in \eqref{Eq:RelThibault} that the
assumptions of Corollary~\ref{Bsub-th} are satisfied when all the
matrices from the generalized Jacobian $\partial_C H(\ox)$ are
nonsingular. In the latter case the solvability of subproblem
\eqref{eq:B-diffnewt} and the superlinear convergence of the
$B$-differentiable Newton method follow from the results of
\cite{QiS 93} that in turn improve the original ones in \cite{Pan
90}, where $H$ is assumed to be strongly Fr\'echet differentiable
at the solution point.

Further, it is shown in \cite{Qi 93} that the $B$-differentiable
method for semismooth equations \eqref{eq:ne} converges
superlinearly to the solution $\ox$ if just matrices
$A\in\partial_B H(\ox)$ are nonsingular while {\em assuming} in
addition that subproblem \eqref{eq:B-diffnewt} is {\em solvable}.
As illustrated by the example presented on pp.\ 243--244 of
\cite{Qi 93}, without the latter assumption the $B$-differentiable
Newton method may not be well defined for semismooth mappings $H$
on the plane with all the nonsingular matrices from $\partial_B
H(\ox)$. We want to emphasize that the solvability assumption for
\eqref{eq:B-diffnewt} is not imposed in
Theorem~\ref{thm:Bdiffconv}---it is {\em ensured} by {\em metric
regularity}.\vspace*{0.05in}

Let us now discuss interconnections between the metric regularity
property of locally Lipschitzian mappings $H\colon\R^n\to\R^n$ via
its coderivative characterization \eqref{cod-cr} and the
nonsingularity of the generalized Jacobian and $B$-subdifferential
of $H$ at the reference point. To this end, observe the following
relationships between the corresponding constructions.

\begin{Proposition}\label{relat} {\bf (relationships between the
$B$-subdifferential, generalized Jacobian, and coderivative of
Lipschitzian mappings).} Let $H\colon\R^n\to\R^m$ be locally
Lipschitzian around $\ox$. Then we have
\begin{equation}\label{inc}
\partial_B H(\ox)^Tz\subset\partial\la z,H\ra(\ox)\subset\partial_C
H(\ox)^Tz\;\mbox{ for all }\;z\in\R^m,
\end{equation}
where both inclusions in \eqref{inc} are generally strict.
\end{Proposition}

\noindent
{\bf Proof.} Recall that the middle term in \eqref{inc} expressed
via the limiting subdifferential \eqref{sub} is exactly the
coderivative $D^*H(\ox)(z)$ due to the scalarization formula
\eqref{Eq:scalarization} for locally Lipschitzian mappings. Thus
the second inclusion in \eqref{inc} follows immediately from the
well-known equality \eqref{C2} involving convexification, and it
is strict as a rule due to the usual nonconvexity of the limiting
subdifferential; see \cite{mor06,rw}.

To justify the first inclusion in \eqref{inc}, observe that the
limiting subdifferential $\partial f(\ox)$ of every function
$f\colon\R^n\to\R$ continuous around $\ox$ admits the
representation
\begin{equation}\label{F}
\partial f(\ox)=\Limsup_{x\to\ox}\Hat\partial f(x)
\end{equation}
via the outer limit \eqref{pk} of the Fr\'echet/regular
subdifferentials
\begin{equation}\label{F1}
\Hat\partial f(x):=\Big\{p\in\R^n\Big|\;\liminf_{u\to
x}\frac{f(u)-f(x)-\la p,u-x\ra}{\|u-x\|}\ge 0\Big\}
\end{equation}
of $f$ at $x$; see, e.g., \cite[Theorem~1.89]{mor06}. We obviously
have from \eqref{F1} that $\Hat\partial f(\ox)=\{f'(\ox)\}$ if $f$
is (Fr\'echet) differentiable at $\ox$ with its
derivative/gradient $f'(\ox)$.

Having the mapping $H=(h_1,\ldots,h_m)\colon\R^n\to\R^m$ in the
proposition and fixing an arbitrary vector
$\oz=(\oz_1,\ldots,\oz_m)\in\R^m$, form now a scalar function
$f_{\oz}\colon\R^n\to\R$ by
\begin{equation}\label{f1}
f_{\oz}(x):=\sum_{i=1}^m\oz_i h_i(x),\quad x\in\R^n.
\end{equation}
Then the first inclusion in \eqref{inc} amounts to say that
\begin{equation}\label{inc1}
   \partial_B H(\ox)^T\oz\subset\partial f_{\oz}(\ox).
\end{equation}
To proceed with proving \eqref{inc1}, pick any matrix
$A\in\partial_B H(\ox)^T\oz$ and denote by $a_i\in\R^n$,
$i=1,\ldots,n$, its vector rows. By definition \eqref{eq:B-sub} of
the $B$-subdifferential $\partial_B H(\ox)$ there is a sequence
$\{x^k\}\subset S_H$ from the set of differentiability \eqref{R}
such that $x^k\to\ox$ and $H'(x^k)\to A$ as $k\to\infty$. It is
clear from \eqref{f1} that the function $f_{\oz}$ is
differentiable at each $x^k$ with
\begin{equation*}
f'_{\oz}(x^k)=\sum_{i=1}^m\oz_i h'_i(x^k)\to\sum_{i=1}^m\oz_i
a_i=A^T\oz\;\mbox{ as }\;k\to\infty.
\end{equation*}
Since $\Hat\partial f_{\oz}(x^k)=\{f'_{\oz}(x^k)\}$ at all the
points of differentiability, we arrive at \eqref{inc1} by
representation \eqref{F} of the limiting subdifferential and thus
justify the first inclusion in \eqref{inc}.

To illustrate that the latter inclusion may be strict, consider
the function $H(x):=|x|$ on $\R$. Then $\partial_B H(0)z=\{-z,z\}$
for all $z\in\R$, while
\begin{eqnarray*}
\partial(z H)(0)=D^*H(0)(z)=\left\{\begin{array}{ll}
[-z,z] &\mbox{for }\;z\ge 0,\\
\{-z,z\}&\mbox{for }\;z<0.
\end{array}\right.
\end{eqnarray*}
This completes the proof of the proposition. $\h$\vspace*{0.05in}

It follows from Proposition~\ref{relat} in the case of
Lipschitzian transformations $H\colon\R^n\to\R^n$ that the
nonsingularity of all the matrices $A\in\partial_C H(\ox)$ is a
sufficient condition for the metric regularity of $H$ around $\ox$
due to the coderivative criterion \eqref{cod-cr} while the
nonsingularity of all $A\in\partial_B H(\ox)$ is a necessary
condition for this property. Note however, as it has been
discussed above, that the nonsingularity condition for $\partial_B
H(\ox)$ alone does not ensure the solvability of subproblem
\eqref{eq:B-diffnewt} in the $B$-differentiable Newton method, and
thus it cannot be used alone for the justification of algorithm
\eqref{eq:B-diffnewt}, \eqref{eq:iter} in the $B$-differentiable
semismooth case. Furthermore, we are not familiar with any
verifiable condition to support the nonsingularity of $\partial_B
H(\ox)$ in the full justification of the $B$-differentiable Newton
method.

In contrast to this, the metric regularity itself---via its
verifiable pointwise characterization \eqref{cod-cr}---ensures the
solvability of \eqref{eq:B-diffnewt} and fully justifies the
B-differentiable Newton method with its superlinear convergence
provided that the mapping $H$ is semismooth and locally invertible
around the reference solution point. Note that the nonsingularity
of the generalized Jacobian $\partial_C H(\ox)$ implies not only
the metric regularity but simultaneously the semismoothness and
local invertibility of a Lipschitzian transformation
$H\colon\R^n\to\R^n$. However, the latter condition fails to spot
a number of important situations when all the assumptions of
Theorem~\ref{thm:Bdiffconv} are satisfied; see, in particular,
Corollary~\ref{Bsub-th} and the corresponding conditions in terms
of Warga's derivate containers discussed right after
Corollary~\ref{suffH1}. We refer the reader to the specific
mappings $H\colon\R^2\to\R^2$ from \cite[Example~2.2]{k} and
\cite[Example~3.3]{w} that can be used to illustrate the above
statement.

\section{Concluding Remarks}

In this paper we develop a new generalized Newton method for
solving systems of nonsmooth equations $H(x)=0$ with
$H\colon\R^n\to\R^n$. Local superlinear convergence and global (of
the Kantorovich type) convergence results are derived under
relatively mild conditions. In particular, the local Lipschitz
continuity and directional differentiability of $H$ are not
necessarily required. We show that the new method and its
specifications have some advantages in comparison with previously
known results on the semismooth and $B$-differentiable versions of
the generalized Newton method for nonsmooth Lipschitz equations.

Our approach is heavily based on advanced tools of variational
analysis and generalized differentiation. The algorithm itself is
built by using the graphical/contingent derivative of $H$, while
other graphical derivatives and coderivatives are employed in
formulating appropriate assumptions and proving solvability and
convergence results. The fundamental property of metric regularity
and its pointwise coderivative characterization play a crucial
role in the justification of the algorithm and its satisfactory
performance.

In the other lines of developments, it seems appealing to develop
an alternative Newton-type algorithm, which is constructed by
using the basic coderivative instead of the graphical derivative.
This requires certain symmetry assumptions for the given problem,
since the coderivative is an extension of the adjoint derivative
operator. Major advantages of a coderivative-based Newton method
would be comprehensive calculus rules held for the coderivative in
contrast to the contingent derivative, complete coderivative
characterizations of Lipschitzian stability, and explicit
calculations of the coderivative in a number of settings important
for applications. The details of these ideas are part of our
future research.

\end{document}

%% file: Hfun.eepic
\setlength{\unitlength}{0.00087489in}
\begingroup\makeatletter\ifx\SetFigFont\undefined%
\gdef\SetFigFont#1#2#3#4#5{%
  \reset@font\fontsize{#1}{#2pt}%
  \fontfamily{#3}\fontseries{#4}\fontshape{#5}%
  \selectfont}%
\fi\endgroup%
{\renewcommand{\dashlinestretch}{30}
\begin{picture}(9180,8589)(0,-10)
\texture{44555555 55aaaaaa aa555555 55aaaaaa aa555555 55aaaaaa aa555555 55aaaaaa 
	aa555555 55aaaaaa aa555555 55aaaaaa aa555555 55aaaaaa aa555555 55aaaaaa 
	aa555555 55aaaaaa aa555555 55aaaaaa aa555555 55aaaaaa aa555555 55aaaaaa 
	aa555555 55aaaaaa aa555555 55aaaaaa aa555555 55aaaaaa aa555555 55aaaaaa }
\put(8112,8067){\shade\ellipse{128}{128}}
\put(8112,8067){\ellipse{128}{128}}
\put(4512,4512){\shade\ellipse{128}{128}}
\put(4512,4512){\ellipse{128}{128}}
\put(2712,2712){\shade\ellipse{128}{128}}
\put(2712,2712){\ellipse{128}{128}}
\path(12,912)(8562,912)
\blacken\path(8442.000,882.000)(8562.000,912.000)(8442.000,942.000)(8478.000,912.000)(8442.000,882.000)
\path(912,12)(912,8562)
\blacken\path(942.000,8442.000)(912.000,8562.000)(882.000,8442.000)(912.000,8478.000)(942.000,8442.000)
\path(8112,912)(8112,732)
\path(4512,912)(4512,732)
\path(2712,912)(2712,732)
\dashline{60.000}(912,912)(8112,8112)
\path(732,2712)(912,2712)
\path(732,4512)(912,4512)
\path(732,8112)(912,8112)
\dashline{60.000}(8112,6762)(507,1092)
\dashline{60.000}(8112,6312)(417,2487)
\dashline{60.000}(2712,912)(2712,2712)
\dashline{60.000}(4512,912)(4512,4512)
\dashline{60.000}(8112,912)(8112,8112)
\path(8112,8112)(5772,5142)(5052,5052)
	(4782,4647)(4647,4647)
\path(4512,4512)(3477,3297)(3207,3207)
	(3117,3072)(3072,3027)
\put(2592,417){\makebox(0,0)[lb]{{\SetFigFont{12}{14.4}{\rmdefault}{\mddefault}{\updefault}$1/4$}}}
\put(4377,417){\makebox(0,0)[lb]{{\SetFigFont{12}{14.4}{\rmdefault}{\mddefault}{\updefault}$1/2$}}}
\put(8067,417){\makebox(0,0)[lb]{{\SetFigFont{12}{14.4}{\rmdefault}{\mddefault}{\updefault}$1$}}}
\put(427,8067){\makebox(0,0)[lb]{{\SetFigFont{12}{14.4}{\rmdefault}{\mddefault}{\updefault}$1$}}}
\put(307,4407){\makebox(0,0)[lb]{{\SetFigFont{12}{14.4}{\rmdefault}{\mddefault}{\updefault}$1/2$}}}
\put(307,2612){\makebox(0,0)[lb]{{\SetFigFont{12}{14.4}{\rmdefault}{\mddefault}{\updefault}$1/4$}}}
\put(7347,6807){\makebox(0,0)[lb]{{\SetFigFont{12}{14.4}{\rmdefault}{\mddefault}{\updefault}$H$}}}
\put(5862,6177){\makebox(0,0)[lb]{{\SetFigFont{12}{14.4}{\rmdefault}{\mddefault}{\updefault}$x$}}}
\put(8327,6707){\makebox(0,0)[lb]{{\SetFigFont{12}{14.4}{\rmdefault}{\mddefault}{\updefault}$\frac{3}{4}x+\frac{1}{16}$}}}
\put(8327,6212){\makebox(0,0)[lb]{{\SetFigFont{12}{14.4}{\rmdefault}{\mddefault}{\updefault}$\frac{1}{2}x+\frac{1}{4}$}}}
\end{picture}
}